\documentclass[12pt, a4paper]{amsart}
\usepackage{amssymb, amsthm, amsmath, amsfonts}

\usepackage[marginpar=1cm, margin=2.6cm]{geometry}
\usepackage{verbatim}

\theoremstyle{plain}

\begin{document}

\title{A remark on a proof of $[G, L]=0$ for a Lie group $G$}

\author{Haruo Minami}
\address{H. Minami: Professor Emeritus, Nara University of Education}
\email{hminami@camel.plala.or.jp}


\begin{abstract}
In this note we give an improvement of our proof of $[G, L]=0$ 
for a compact framed Lie group $(G, L)$, which depends heavily 
on the choice of a circle subgroup $S\subset G$. We attempt 
here to make a more suitable choice of this circle subgroup.
\end{abstract}

\maketitle

\section{Introduction}

Let $G$ be a compact connected Lie group and $L$ be 
its left invariant framing. We denote by $[G, {\it L}]$ the framed 
bordism class of $(G, L)$. In ~\cite{M} we gave a proof of the 
following result, based on a proof technique proposed in ~\cite{O}.

{\sc Theorem.}  {{\it For} \ $G=SU(n)$, \ 
$SO(n)$, \ $Spin(n)$ \ $(n \ge 8)$; \ $Sp(n)$ \ $(n \ge 4)$; \ 
$F_4,$, \ $E_6$, \ $E_7$, \ $E_8$, we have $[G, L]=0$.}

Here we give a direct proof by making a more 
suitable choice of the circle subgroup $S\subset G$ 
which is a key ingredient of the formula used in ~\cite{O}. 

From \cite{O, K} we know that $[G, L]_{(p)}=0$ for all primes 
$p\ge 5$, so we restrict ourselves to the cases $p=2, 3$ where 
$x_{(p)}$ is the image of $x$ by the localizing 
homomorphism at $p$.

Suppose we are given a circle subgroup $S\subset G$ and a 
one-dimensional representation $\gamma : S \to U(1)$.
Let $\xi$ be the complex line bundle associated with the 
principal $S$-bundle $G \to G/S$ by $\gamma$. Assuming   
that its sphere bundle $S(\xi) \to G/S$ is isomorphic to 
$G\to G/S$, we consider the Kronecker product of $[G/S]\in 
\pi_{d-1}^S(G/S^+)$ and $J(\beta\xi)\in \pi_S^d(S^1(G/S^+))$ 
($d=\dim G$). Then by  ~\cite{O} we have
\[[G, L]=-\langle J(\beta\xi)), [G/S]\rangle \quad \text{in}  
\ \pi_d^S.\]
Here $G/S$ is a framed manifold with the framing inherited   
from $G$ in the natural way and let $\beta$ be the Bott 
element and $J$ be the complex $J$ homomorphism.

In view of this formula, in order to prove that $[G, L]=0$ we show that 
$J(\beta\xi)=0$ holds for $S\subset G$ specified depending on each 
$G$. But in fact for the reasons mentioned above we show that 
$J(\beta\xi)_{(p)}=0$ holds only for $p=2, 3$. 

For $j \in \mathbf{Z}$ we set $t^j=J(\beta\xi^j)$ (where $\xi^0$ 
is the trivial line bundle $G/S\times \mathbf{C}$). Then 
by the solution of the Adams conjecture we have \begin{equation*}
t_{(p)}^j = kt_{(p)}^{kj} \quad (k, j \ne 0) \quad \text{if} \ (p, k)=1. 
\end{equation*} 
Also, since $J(\beta)$ becomes a generator of $\pi^S_1=\mathbf{Z}_2$, 
we have  
\begin{equation*} 
2\!\cdot\!1=0 
\end{equation*}
where $1=t^0$. Applying these two relations to (1) below enables 
us to calculate that $t_{(p)}=0$. But in the calculation below we drop 
for brevity the subscript of $t_{(p)}$ except in some exceptional cases.

For any given $k_1, \cdots, k_\ell\in \mathbf{Z}$ we write   
$S(k_1, \cdots, k_\ell)$ for the circle subgroup of $G$ generated by 
${\rm diag}(z^{k_1}, \cdots, z^{k_\ell})$, $z\in U(1)$, where 
${\rm diag}(c_1,  \cdots, c_\ell)$ denote the diagonal matrix whose 
$ii$-entry is $c_i$. We take $S=S(k_1, \cdots, k_\ell)$ for 
$k_1, \cdots, k_\ell\in \mathbf{Z}$ with $k_j=1$ for some $j$ and  
then define $\gamma : S \to U(1)$ by 
${\rm diag}(z^{k_1}, \cdots, z^{k_\ell})\mapsto z$. Suppose here
that there exits a complex representation $\rho : G \to U(n)$ satisfying    
\begin{equation*}
\rho|S=\gamma^{k_1}+\cdots+\gamma^{k_\ell}+(n-s) \quad  (s\ge 0). 
\end{equation*}
Then it holds that $\xi^{k_1}\oplus\cdots\oplus\xi^{k_\ell}
\oplus (n-s)\xi^0\cong n\xi^0$ and therefore we have
\begin{equation} 
t^{k_1+i}+ \cdots +t^{k_\ell+i}=\ell t^i \quad (i \ge 0) 
\end{equation}
Below we write $(1)_i$ for this equality to clarify that it belongs 
to the formula of $t^i$. Besides in some cases, this equality is used in 
combination with the one obtained by using $\lambda^j\rho$ 
instead of $\rho$. 
.  

\section{Proof of Theorem for classical Lie groups}

\begin{proof}[Proof for the case $G=Sp(n)$]  Let $Sp(n)$ be 
embedded in $SU(2n)$ in the standard way. Let $\rho :Sp(n) 
\to U(2n)$ be the restriction of the inclusion homomorphism  
$SU(2n) \to U(2n)$ to $Sp(n)$. Take 
$S=S(1, 2, 3, -6, 0, \cdots, 0)\subset Sp(n)$. Then this circle 
subgroup corresponds to 
$S(1, -1, 2, -2, 3, - 3, -6,  6, 0, \cdots, 0)$ in 
$SU(2n)$ via the embedding above, so by (1) we have
\begin{equation} 
t^{1+i}+t^{2+i}+t^{3+i}+t^{6+i}+t^{-1+i}+t^{-2+i}+t^{-3+i}+t^{-6+i}=8t^i \quad 
(i \ge 0). 
\end{equation}  
{\it Case} $p=2$. From ${(2)}_1$, ${(2)}_2$, ${(2)}_4$ we have
\[
105t^4=916t+1, \ 45t^4+60t^2=76t+1, \ 840t^4-56t^2=176t  
\]
(where the subscript of $t_{(2)}$ is omitted as noted above). By 
eliminating $t^2$ and $t^4$ from these equations we have \[16t=0\]  
and therefore $t^4=4t+1$, $16t^2=0$. Substituting these 
equalities into ${(2)}_5$, ${(2)}_2$ we have $4t^2+8t=0$, $t^8=8t+1$. 
Finally, substituting all these equalities into $(2)_6$ we obtain  
$t=0$, i.e., $t_{(2)}=0$.

{\it Case} $p=3$. From $(2)_1$, $(2)_2$, $(2)_3$ (with $t$ 
replaced by $t_{(3)}$) we have similarly 
\[3t=0, \ t^3=0, \ t^9=0.\] In the above, replacing $\rho$ by 
$\lambda^2\rho$, we have a similar equality to (2):
\begin{equation*}
\begin{array}{l}   
2t^{1+i}+t^{2+i}+2t^{3+i}+2t^{4+i}+2t^{5+i}+t^{7+i}+t^{8+i}+t^{9+i}
+2t^{-1+i}\\+t^{-2+i}+2t^{-3+i}+2t^{-4+i}+2t^{-5+i}+t^{-7+i}+t^{-8+i}
+t^{-9+i}=24t^i \quad (i \ge 0). 
\end{array}
\end{equation*}
By substituting the three equalities obtained above into this equality 
for $i=1$ we obtain $t=0$, i.e., $t_{(3)}=0$.  
\end{proof}

\begin{proof}[Proof for the case $G=SU(n)$]   Let $\rho$ be 
the inclusion homomorphism $SU(n) \to U(n)$. We take  
$S=S(1, -1, 2, -2, 3, - 3, -6,  6, 0, \cdots, 0)\subset SU(n)$ which  
has the same form as that in the case $G=Sp(n)$. Then  
this choice of $S$ allows us to apply the same argument and 
consequently leads us to the result.     
\end{proof}

\begin{proof}[Proof for the case $G=SO(n)$]  Let $\rho : SO(n) \to 
U(n)$ be the complexification of the real inclusion homomorphism 
$SO(n) \to O(n)$. Let $S\subset SO(n)$ be the circle subgroup 
corresponding to $S(1, -1, 2, -2, 3, - 3, -6,  6, 0, \cdots, 0)$ in 
$SU(n)$ via the canonical embedding $\iota : SO(n) \to SU(n)$. 
For the same reason as for the above, applying the argument as 
in the above case we can obtain the result similarly.
\end{proof}

\begin{proof}[Proof for the case $G=Spin(n)$]  Let $\rho=i\circ\pi : 
Spin(n) \to U(n)$ where $\pi$ denotes the natural covering morphism 
$Spin(n) \to SO(n)$ and $i : SO(n) \to U(n)$ the complexification of 
the inclusion homomorphism $SO(n) \to O(n)$. Let $S\subset 
Spin(n)$ be the circle subgroup such that its image by $\pi$ 
corresponds to $S(2, -2, 2, -2, 4, -4, 6, -6, 0, \cdots, 0)$ in $SU(n)$ 
via $\iota$ above. It is clear that by definition $S$ contains 
$-1\in Spin(n)$. Therefore by applying (1) we have 
\begin{equation} 
2t^{2+i}+t^{4+i}+t^{6+i}+2t^{-2+i}+t^{-4+i}+t^{-6+i}=8t^i \quad 
(i \ge 0). 
\end{equation}

{\it Case} $p=2$.  From $(3)_1$ we have 
\[8t=0.\] 
In the above, replacing $\rho$ by $\lambda^4\rho$ we have 
\begin{equation*}
\begin{array}{l}   
t^{2+i}+2t^{4+i}+5t^{6+i}+4t^{8+i}+4t^{10+i}+t^{14+i}+t^{-2+i}\\+2t^{-4+i}
+5t^{-6+i}+4t^{-8+i}+4t^{-10+i}+t^{-14+i}=34t^i \quad (i \ge 0). 
\end{array}
\end{equation*} 
Substituting $8t=0$ above into this equality for $i=1$ 
we have $t=0$, i.e., $t_{(2)}=0$. 

{\it Case} $p=3$. From $(3)_1$ we have $t_{(3)}=0$ by a simple calculation.
\end{proof}

\section{Proof of Theorem for exceptional Lie groups}

\begin{proof}[Proof for the case $G=F_4$] 
We know ~\cite{O} that $F_4$ has $Spin(9)$ as a subgroup and   
a representation $U : F_4 \to U(26)$ such that its restriction to this 
$Spin(9)$ is $1+ \lambda^1+\Delta$, where $\lambda^1=i\circ\pi : 
Spin(9) \to SO(9) \to U(9)$ (with the notation above) and $\Delta$ is the 
spin representation. Take $\rho=U$. Then, if we choose $S\subset F_4$ 
so that its image by $\pi$ corresponds to $S(2. -2, 2. -2, 2. -2, 4. -4, 0)$ 
in $SU(9)$ via the canonical embedding $\iota$, then we have 
\[\rho\!\mid\!Spin(9)=(\gamma+\gamma^{-1})^3(\gamma^2
+\gamma^{-2})+3\gamma^2+3\gamma^{-2}+\gamma^4
+\gamma^{-4},\] so by (1) we have 
\begin{equation}
\begin{array}{l}
4t^{1+i}+3t^{2+i}+3t^{3+i}+t^{4+i}+t^{5+i}+4t^{-1+i}+3t^{-2+i}\\+3t^{-3+i}
+t^{-4+i}+t^{-5+i}=24t^i \quad (i \ge 0). 
\end{array} 
\end{equation}
By replacing $\rho$ by $\lambda^2\rho$ we also have
\begin{equation*}
\tag{4'}
\begin{array}{l}   
25t^{1+i}+24t^{2+i}+19t^{3+i}+19t^{4+i}+13t^{5+i}+10t^{6+i}+6t^{7+i}
+3t^{8+i}\\+t^{9+i}+25t^{-1+i}+24t^{-2+i}+19t^{-3+i}+19t^{-4+i}
+13t^{-5+i}+10t^{-6+i}\\+6t^{-7+i}+3t^{-8+i}+t^{-9+i}=240t^i \quad (i \ge 0). 
\end{array}
\end{equation*}

{\it Case} $p=2$. From $(4)_1$, $(4)_2$ we have 
$20t^2+30t^4=392t$, $2590t^2-315t^4=288t$, respectively. 
Further, calculating both $(4)_3$ and $(4)_4$ we have $7840t^4=3128t$. 
By eliminating $t^2$, $t^4$ from these equations we have $8t=0$, so 
$32t^2=0$, $2t^2=5t^4$. Calculating $(4)_3$ and $(4)_4$ again 
by use of these equalities, we have  
\[8t=0, \ 8t^2=0, \ t^4=2t^2, \ t^8=4t^2+1.\]
Using $(4')_1$ the second equality is refined into $4t^2=0$, so 
it follows that $ t^8=1$. Substituting these equalities
consequently we have $t=0$, i.e., $t_{(2)}=0$.

{\it Case} $p=3$. From $(4)_1$, $(4)_2$ we have $3t=0$, $t^3=0$. 
Using these equalities, from $(4)_4$ we have $t^9=0$. Substituting 
these equalities into $(4')_1$ yields $t=0$, i.e., $t_{(3)}=0$.
\end{proof}

\begin{proof}[Proof for the case $G=E_6$] 
From ~\cite{O} we know that $E_6$ has a subgroup $F_4$ 
and a representation $W : E_6 \to U(27)$ such that its restriction to 
$F_4\subset E_6$ is $1+U$. 
This means that it enables us to apply the proof of the case $G=F_4$ 
to the case here. Consequently we obtain the result.
\end{proof}

\begin{proof}[Proof for the case $G=E_7$] By ~\cite[Theorem. 11.1]{O}, 
$E_7$ contains $SU(8)/\{\pm I\}$ as a sybgroup where $I\in SU(8)$ is 
the identity and a representation $\rho  : E_7 \to U(56)$ such 
that its restriction to this subgroup group is $\lambda^2+\lambda^4$. 
Here $\lambda^j$ denotes the $j$-th exterior power of the standard 
representation of $SU(8)$ on $\mathbf{C}^8$. Take 
$S=S(1, \cdots 1,-7)/\{\pm I\} \subset E_7$. 
Then by (1) we have
\begin{equation}
21t^{1+i}+7t^{3+i}+21t^{-1+i}+7t^{-3+i}=52t^i\quad (i\ge 0).
\end{equation}
Here $t$ is the replacement of the square of $t$ defined for 
$S(1, \cdots, 1-7)$ in (1). But by definition of (1) we find that in order to obtain 
the required result it suffices to prove that $t_{(2)}=0$ and $t_{(3)}=0$ for 
this $t$. 

In addition, replacing $\rho$ by $\lambda^3\rho$ in the above, 
we also have  
\begin{equation*}
\tag{5'}
\begin{array}{l}   
3t^{1+i}+t^{3+i}+6t^{5+i}+t^{7+i}+3t^{9+i}+6t^{-1+i}
+3t^{-1+i}\\+t^{-3+i}+6t^{-5+i}+t^{-7+i}+3t^{-9+i}+6t^{1+i}=0 \ \mod 8 
\qquad (i\ge 0).
\end{array}
\end{equation*}

{\it Case} $p=2$. From $(5)_1$, $(5)_2$ and $(5)_3$ we have $16t=0$, so 
$8t^2=0$ and $t^4+2t^2=8t+1$. Using $(5)_4$ the first equality is 
refined into $4t=0$, so the last one becomes $t^4+2t^2=1$. From $(5)_5$ 
we have $t^8=4t^2+1$. Substituting these equalities, from $(5')_3$, $(5')_4$ 
we have $1=0$ and $t=1$, respectively. Combining these two results 
we have $t=0$, i.e,. $t_{(2)}=0$. 

{\it Case} $p=3$. From the calculation of $(5)_1$, $(5)_2$, $(5)_3$ we have 
$t=0$, i.e., $t_{(3)}=0$. 
\end{proof}

\begin{proof}[Proof for the case $G=E_8$] We know ~\cite{A} that 
$E_8$ contains $Spin(16)$ as a subgroup and the restriction of 
the adjoint representation of $E_8$ to this $Spin(16)$ is 
$\lambda^2+\Delta^+$ where $\lambda^2$ is the adjoint 
representation of $Spin(16)$ and $\Delta^+$ the 
positive spinor representation. We choose here a different 
$S\subset E_8$ in each case.

{\it Case} $p=2$. Let $S\subset E_8$ be the circle subgroup of 
$Spin(16)$ such that its image by $\iota\circ\pi : Spin(16) \to SO(16) 
\to SU(16)$ is $S(2, -2, \cdots, 2, -2, 6, -6, 0, 0)\subset 
SU(16)$. Then by (1) we have 
\begin{equation}
\begin{array}{l}
21t^{1+i}+12t^{2+i}+21t^{3+i}+21t^{4+i}+15t^{5+i}+2t^{6+i}+6t^{7+i}
+6t^{8+i}+t^{9+i}\\+21t^{-1+i}+12t^{-2+i}+21t^{-3+i}+21t^{-4+i}
+15t^{-5+i}+2t^{-6+i}+6t^{-7+i}
\\+6t^{-8+i}+t^{-9+i}=210t^i \quad (i \ge 0).    
\end{array}  
\end{equation}
Calculating $(6)_1$, $(6)_2$, $(6)_3$, $(6)_4$ we have 
\[2^5\cdot 10223201t^2=2\cdot 1331664213t, \  \ 
2^2\cdot 248349949873t^2=2^2\cdot 104608953537t.\] 
From these equalities it follows that $2t=0$ and thereby 
$4t^2=0$, $2t^4=0$, $t^8=0$. Further, by using these equalities, 
from $(6)_7$ we have $t^{16}=2t^2+1$. Finally, substituting these 
equalities into $(6)_8$ we obtain $t=0$, i.e., $t_{(2)}=0$.

{\it Case} $p=3$. We choose the circle subgroup of $Spin(16)$ as 
$S\subset E_8$ such that its image by $\iota\circ\pi : Spin(16) 
\to SO(16) \to SU(16)$ is $S(2, -2, \cdots, 2, -2, -2, 2)\subset 
SU(16)$. Then by (1) we have
\begin{equation}
56t^{2+2i}+28t^{8+2i}+8t^{6+2i}+56t^{-2+2i}+28t^{-8+2i}+8t^{-6+2i}
=184t^{2i} \quad (i \ge 0).    
\end{equation} 
Here, thinking of $t^2$ as $t$ we prove that $t_{(3)}^2=0$, which 
means that $t_{(3)}=0$, because of $2t_{(3)}^2=t_{(3)}$,

Combining the resuts of calculations of ${(7)}_{2i}$ for $i=1, 2, 3, 4$, 
we have $3t^2=0$, $t^6=0$. Calculating ${(7)}_{10}$ by use of  
these equalities we have $t^{18}=0$. Substituting all these 
equalities into ${(7)}_{14}$, we obtain $t^2=0$, i.e., $t^2_{(3)}=0$, so 
as stated above we conclude that $t_{(3)}=0$.   
This completes the proof of the theorem.
\end{proof}

\end{document}